

\input amstex
\documentstyle{amsppt}

\refstyle {A}
\widestnumber\key{AAA}
\magnification=1200
\baselineskip=15pt
\parskip 3pt
\pagewidth{5.5in}
\pageheight{7.4in}

\def\today{\ifcase\month\or
  January\or February\or March\or April\or May\or June\or
  July\or August\or September\or October\or November\or
  December\fi
  \space\number\day, \number\year}

\def\leq{\leqslant}
\def\geq{\geqslant}
\def\eps{\varepsilon}
\def\OO{{\Cal O}}
\def\Ab{{\frak A\!b}}
\def\NS{\text{\rm NS}}
\def\RR{{\Bbb R}}
\def\PP{{\Bbb P}}
\def\phi{\varphi}

\topmatter
\title Remarks on Seshadri constants  \endtitle
\author Andreas Steffens \endauthor
\address{ Andreas Steffens \hfill {} \linebreak
\hglue 13.6pt Department of Mathematics \hfill {} \linebreak
\hglue 13.6pt University of California \hfill {} \linebreak
\hglue 13.6pt Los Angeles, CA 90024}  \endaddress
\email  steffens\@math.ucla.edu \endemail
\date \today \enddate
\thanks Supported by Deutsche Forschungsgemeinschaft \endthanks
\abstract
Given a smooth complex projective variety $X$ and an ample line bundle $L$ on
$X$.
Fix a point $x\in X$.
We consider the question, are there conditions which guarantee the maxima of
the Seshadri constant of $L$ at $x$, i.e  $\eps(L,x)=\root n \of {L^n}$?
We give a partial answer for surfaces and find  examples where the answer to
our question
is negative. If $(X,\Theta)$ is a general principal polarized abelian surface,
then $\eps(\Theta,x)=\frac{4}{3}<\sqrt{2}=\sqrt{\Theta^2}$ for all $x\in X$.
\endabstract
\endtopmatter

\document

\heading  Introduction \endheading

Let $X$ be a smooth projective variety and let $L$ be a line bundle on $X$.
Fix a point $x\in X$.
Demailly \cite{D} introduced a very interesting measure of the local positivity
at a point $x$ of  $L$, namely the real number
$$\eps(L,x)= {\mathrel{\mathop{\text{\rm inf}}\limits_{C\ni x}}}\
\frac{L.C}{m_x(C)},$$
which is called the {\bf Seshadri constant} of $L$ at $x$.
Here the infimum is taken over
all irreducible curves $C$ passing through $x$ and $m_x(C)$ is the
multiplicity of $C$ at $x$.
For example, if $L$ is very ample then $\eps(L,x)\ge 1$.

There has been recant interest in trying to give lower bounds for this
invariant at a general point.
Ein and Lazarsfeld \cite{EL} show that if $X$ is a surface, then $\eps(L,x)\ge
1$
for very general $x\in X$.
 In higher dimension
$(n \geq 3)$ Ein, K\"uchle and Lazarsfeld \cite{EKL} prove that
$\eps(L,x)\geq\frac{1}{n}$
for a very general point.
We say that a point $x\in X$ is very general if $x$ is in the complement
$X \smallsetminus Z$ of $Z$  a
countably union of proper subvarieties.
 Examples of Miranda show that
$\eps(L,x)$ can take arbitrarily small values in codimension two
(i.e. $\text{codim}\ Z=2$),
even for an  ample line bundle.

One may expect that this general bounds are not optimal. An elementary
observation
(see Remark 1 below for the proof) shows that $\eps(L,x)\leq \root n \of
{L^n}$.
A natural question is,
are there conditions which guarantee  equality?
Even in relative simple  cases it turns out to be hard to give
an answer.

Recently, Xu \cite{Xu} improved the surface bound given by Ein and Lazarsfeld.
He showed that if $L^2\geq \frac{1}{3}(4 \alpha^2-4\alpha + 5)$
for a given integer $\alpha >1$
and $L.C\geq \alpha$ for every irreducible curve $C\subset X$,
then $\eps(L,x)\geq \alpha$ for all but finitely many $x\in X$.

The first result we have gives a further improvement and gives a partial answer
to the equality question.

\proclaim{Proposition 1} Let $X$ be a surface with $\rho(X)=rk(\NS(X))=1$
and let $L$ be an
ample generator of $\NS(X)$. Let $\alpha$ be an integer with $\alpha^2 \leq
L^2$.

If $x \in X$ is a very general point, then $\eps(L,x)\geq \alpha$.
In particular if $\sqrt{L^2}$ is an integer, then $\eps(L,x)=\sqrt{L^2}.$
\endproclaim

For example, if $(X,L)$ is a general polarized abelian surface of type
$(1,2d^2)$ for some $d\ge 1$,
then $\eps(L,x)=2d=\sqrt{L^2}$ for all $x\in X$.
Or if $X\subset\PP^3$ is a general hypersurface of degree
$d^2\ge 4$, then $\eps(\OO_X(1),x)= d$ for a very general $x\in X$.

The proof uses essentially the fact that $L$ is the  generator and in a
diophantine way that $\alpha$ is an integer.

One might be tempted to suppose that the conclusion of  Proposition 1 holds
allowing $\alpha$ to be a (possibly non-integral) real number.
But the next result shows that the situation is  more complicated.

\proclaim{Proposition 2}
 Let $X$ be the Jacobian of a hyperelliptic curve of genus $g\geq 2$ with
$rk(\NS(X))=1$.
And let $\Theta$ be  the theta divisor on $X$.

Then $\eps(\Theta)=\eps(\Theta,x)\leq \frac{2g}{g+1} <
{\root g\of {g!}}=\root g \of {\Theta^g}$.

\noindent
In particular if $X$ is an irreducible principal polarized abelian surface,
then
$$\eps(\Theta)=\eps(\Theta,x)=\frac{4}{3} <
{\sqrt{2}}=\sqrt{\Theta^2}$$
\endproclaim

Proposition 2 gives an example that the bound
$\eps(L,x)\geq [\sqrt{L^2}] $ does $\underline{not}$ hold for an arbitrarily
line bundle $L$,
 where $[r]$ denotes the integer part of an
real number $r$.
In fact, if $X$ is the Jacobian of a very general curve of genus two, than
 $rk(\NS(X))=1$.
And if $L=\nu \Theta$ then
$\eps(L,x)=\nu\eps(\Theta,x)\leq \nu\frac{4}{3} < [\nu\sqrt{2}]$
for all $\nu\geq 8$.

The arguments in the proof of Proposition 2 works more general to show:

\proclaim{Proposition 3}
Let $X$  be a general principal polarized abelian variety of dimension $g$ with
theta divisor $\Theta$.
Then
$$\eps(\Theta)=\eps(\Theta,x)\le
\root {g-1} \of {\frac{g!}{2^{g-3}(2^g-1)}} <
{\root g \of {g!}} = {\root g \of {\Theta^g}}.$$
\endproclaim

In the case of surface, the following result shows that if $\eps(L,x)$ is
non-maximal (i.e. $\eps(L,x) < \sqrt{L^2}$) then it is rational. There seems to
be no examples known where $\eps(L,x)$ is irrational.

\proclaim{Proposition 4}
In dimension $n$, if the Seshadri constant is non-maximal, then it is a $d$-th
root
of a rational number for some $1\leq d \leq n-1$.
\endproclaim

\subhead Acknowledgment \endsubhead
The author is grateful to R. Lazarsfeld for his great hospitality,
warm encouragements and introducing me to interesting mathematical
question, helpful suggestions and very useful discussions.

\heading 1. Proofs and further remarks \endheading

\demo{Proof of Proposition 1}
Suppose to the contrary that there exists a reduced
irreducible curve $C\subset X$ through a general point $x\in X$, such that
$C.L > \alpha m_x(C)$. Then by the arguments of \cite{EL} it follows that:
$$C^2 \geq m_x(C)(m_x(C)-1). \tag *$$

To see this we follow \cite{EL}. We may assume that $(C,x)$ moves in a
non-trivial
continuous family $\{\ C_t \ni x_t\ \}_{t\in\Delta}$ of reduced irreducible
curves
$C_t\subset X$, plus points $x_t\in C_t$, with
$$m_t=m_{x_t}(C_t) > \alpha C_t.L\ \text{{\hglue 15pt} for all } t\in \Delta.$$
The precise statement we need is:

\proclaim{Propostion\cite{EL}}
Let $\{\ C_t \ni x_t\ \}_{t\in\Delta}$ be a $1$-parameter family of reduced
irreducible curves on a smooth projective surface $X$, such that
$m_t=m_{x_t}(C_t) \geq m$ for all $t\in\Delta$.

Then
$$(C_t)^2 \geq m(m-1).$$
\endproclaim

Now back to the proof of Proposition 1.
By the condition $\rho(X)=1$ there exist an integer $d$, such that $C$ is
numerically equivalent to $dL$. Since $L^2\geq \alpha^2$,  the  assumption
that $C.L > \alpha m_x(C)$ gives
$$\alpha d < m_x(C).\tag **$$
  So it follow from the fact that $\alpha$
is an integer that
$$\alpha d \leq m_x(C)-1. \tag ***$$
Hence by (*),(**) and (***)
$$m_x(C)(m_x(C)-1)\mathrel{\mathop\leq \limits_{(*)}} C^2 = C.(dL)
\mathrel{\mathop < \limits_{(**)}}
 \alpha d m_x(C)
\mathrel{\mathop\leq \limits_{(***)}} m_x(C)(m_x(C)-1),$$
a contradiction.
 \hfill\qed
\enddemo

\demo{Proof of Proposition 2}
Let $C$ be a hyperelliptic curve of genus $g\geq 2$. Then $C$ has $2g+2$
Weierstrass points $p_1, \dots,p_{2g+2}$, with
$2p_1\sim 2p_2\sim \dots \sim 2p_{2g+2}$.  In terms of the Jacobian
$(X,\Theta)=(J(C),\Theta_C)$ of $C$ this has the following interpretation:

Let $\Ab_C: C \rightarrow X\simeq Pic^1(C)$ be the Abel-Jacobi map
and $2_X:X \rightarrow X$ the multiplication by two, the map determine by
the  map
$X\simeq Pic^1(C) \rightarrow X\simeq Pic^2(C), \eta\to cl(\eta^{\otimes 2})$.

Then $C^\prime = 2_X(\Ab_C(C))$ has a point $x$ with $m_x(C^\prime)=2g+2$. And
$C \to C^\prime$ is birational.

 Assume in contrary that the map $f=2_X\Ab_C:C\to C^\prime$ is not birational
and say
has degree $n\geq 2$. Let $\nu :\tilde{C}\to C^\prime$ the normalization of
$C^\prime$.
Then $f$ factors through $\nu$. By the universal property of the Jacobian,
there is a map $\tilde{f} :X=J(C)\to J(\tilde{C})$. Since $rk(\NS(X)=1$  it
follows
that $g=\dim X = \dim J(\tilde{C}) =g(\tilde{C})$.
Hence we find by the Riemann-Hurwitz formula
$$2g-2\geq n(2g(\tilde{C})-2) > 2g-2,$$
a contradiction.

Since $2_X^*\Theta\mathrel{\mathop \sim \limits_{\text{alg}}} 4\Theta$
\cite{LB, II-3 Proposition 3.6}, we find that
$$(\Theta.C^\prime)=(2_X^*\Theta.\Ab_C(C))=
 (4\Theta.\Ab_C(C))=4 g .$$
Hence, we get
$$\eps(\Theta,x)\leq \frac{2g}{g+1} < \root g \of {g!} = \root g \of
{\Theta^g}.$$

Now let  $g=2$ and suppose to the contrary that $\eps(\Theta,x)< \frac{4}{3}$.
Then there exists a reduced irreducible curve $\tilde{C}$ with
$\tilde{C} \mathrel{\mathop{=}\limits_{\text{num}}} aC$ and
$m_x(\tilde{C})=b$ such that
$$\frac{a}{b} < \frac{2}{3}. \tag *$$
Let $\phi : B\!l_x(X) \to X$ the blow-up of $X$ at $x$ with exceptional
divisor $E$. Since $C^\prime$ and $\tilde{C}$ have no common components, it
follows
$$0\leq (\phi^*C^\prime - 6E).(\phi^*\tilde{C}-bE)=8a-6b. \tag **$$
Combining (*) and (**) we find $9a<6b\leq 8a$ the desired contradiction.
\hfill\qed
\enddemo

Sometimes it is useful to use an alternative definition of the Seshadri
constant
of $L$ at a point $x\in X$. If $\phi:B\!l_x(X) \rightarrow X$ is the blow-up
of $X$ at $x$ with exceptional divisor $E$, then
$$\eps(L,x)=\text{sup}\{\ \delta \in \RR \ | \  \phi^*L - \delta E \text{ is
nef }\}$$

\remark{Remark 1}
Let $Y\subset B\!l_x(X)$ be a subvariety of dimension $s=\dim Y$ and
$\delta\leq\eps(L,x)$. Then by Kleiman's theorem \cite{Kl} we have
$(\phi^*L-\delta E)^s.Y\geq 0$. In particular, $(\phi^*L-\delta E)^n\geq 0$.
Hence it follows that $\eps(L,x)\le \root n \of {L^n}$.
\endremark

Let us recall the Nakai-Moishezon criterion for ampleness, which was extended
to
the case of real divisors by Campana and Peternell.
We say that a $\Bbb R$-divisor is ample if its corresponding real point
in the N\'{e}ron-Severi space $N^1(X)$ lies in the interior of the ample cone
of $X$.

\proclaim{Nakai-Moishezon criterion for $\Bbb R$-Cartier divisors \cite{CP}}
\par\noindent
Let $D=\sum a_iD_i$ be a $\Bbb R$-Cartier divisor on a variety $X$.

Then $D$ is ample if and only if $D^s.Y >0$ for any s-dimensional subvariety
$Y\subset X$.
In particular if $D$ is numerically effective but not ample, then there exist
an
irreducible subvariety $Y\subset X$, say of dimension s, such that $D^s.Y=0$.
\endproclaim

\demo{Proof of Proposition 4}
Let $\delta=\eps(L,x)< \root n \of {L^n}$. Then $\phi^*L-\delta E$ is
numerically
effective, but not ample.
Hence by the real Nakai-Moishezon criterion, there exist a subvariety
$Y\subset B\!l_x(X)$
with $(\phi^*L-\delta E)^d.Y=0$, $d=\dim Y$. Since
$(\phi^*L -\delta E)^n > 0$ it follows that $1\leq d \leq n-1$.
Finally by noting that all mixed terms $\phi^*L^i.E^{d-i}$ are zero for
$1\leq i \leq d-1$ we find that $\delta^d$ is rational number.
\hfill\qed
\enddemo

\demo{Proof of Proposition 3}
The general ideology behind the proof is as follows. It might by be difficult
to
bound $\eps(L,x)$ only by using curves, because singular curves are invisible.
Nevertheless, any subvariety $Y\subset X$ with  high multiplicity $m_x(Y)$
at $x$
forces $\eps(L,x)$ to be small.
The precise statement is (c.f.\cite{De, Remark 6.7}):
If $Y$ is a $p$-dimensional subvariety of $X$ passing through $x$ then
$L^p.Y\geq \eps(L,x)^p m_x(Y)$.

Let $X$ be a principal polarized abelian variety of dimension $g$ and let $2_X$
be the multiplication by two.
On $X$ there are $2^{g-1}(2^g -1)$ odd theta characteristics such that
$\Theta$ passe through $2^{g-1}(2^g -1)$ two torsion points
(\cite{Mu, Corollary 3.15 in Appendix to II-3}).
So there is a  divisor
$\Theta^\prime=2_X(\Theta)$  with  a point $x$ having multiplicity
$m=m_x(\Theta^\prime)\geq 2^{g-1}(2^g-1)$ at $x$. And note for late
use that $\Theta^\prime$ is numerically equivalent to $4\Theta$, since
$2_X^*\Theta\mathrel{\mathop \sim \limits_{\text{alg}}} 4\Theta$
\cite{LB, II-3 Proposition 3.6}.

\proclaim{Claim}
$2_X$ maps  $\Theta$ generically $1:1$ to its image.
\endproclaim
Assume to the contrary that  multiplication by two is not generically $1:1$
over $\Theta^\prime$. Then
for general $x\in\Theta$, there is a $y=y(x)\ne x$ such that $2(x-y)=0$.
Then there is a two torsion point $\eta\in X$ such that $(x-y)=\eta$ for all
$x\in X$. But then $\Theta -\eta = \Theta$. But a theta divisor is not
invariant
under any translations.

Now we are in position to compute an upper bound for $\eps(\Theta,x)$, using
the
notation before Remark 1.
Put $\eps=\eps(\Theta,x)$ and let $\widehat\Theta=\phi^*\Theta^\prime-m E$
be the strict transform of $\Theta^\prime$ on $B\!l_x(X)$.
Then by the remark at the beginning we find:
$$0\leq (\phi^*\Theta-\eps E)^{g-1}.\widehat\Theta =
(\phi^*\Theta -\eps E)^{g-1}.(4\phi^*\Theta-m E)= 4(g!) - \eps^{g-1} m.$$
Hence
$$\eps \leq \root {g-1} \of {\frac{4(g!)}{2^{g-1}(2^g-1)}}
= \root {g-1} \of {\frac{g!}{2^{g-3}(2^g-1)}}
  < \root g \of {g!}.$$
\hfill\qed
\enddemo
\bigskip


\enddocument
\end